\documentclass[twocolumn,10pt,final,conference]{IEEEtran}
\usepackage{graphicx}
\usepackage{siunitx}
\graphicspath{{./images/}}
\usepackage{amssymb,amsmath}
\newtheorem{theorem}{Theorem}
\newtheorem{corollary}{Corollary}
\newtheorem{assumption}{Assumption}
\newtheorem{proofoftheorem}{Proof of Theorem}
\title{Feedback Regularization and Geometric PID Control for Trajectory Tracking of Coupled Mechanical Systems: Hoop Robots on an Inclined Plane}

\author{\IEEEauthorblockN{T.\,W.\,U.\,Madhushani\IEEEauthorrefmark{1}, 
D.\,H.\,S.\,Maithripala\IEEEauthorrefmark{2}, 
J.\,M.\,Berg\IEEEauthorrefmark{3}}

\IEEEauthorblockA{\IEEEauthorrefmark{1}Postgraduate and Research Unit, Sri Lanka Technological Campus, CO 10500, Sri Lanka.\\ 
}

\IEEEauthorblockA{\IEEEauthorrefmark{2}Department of Mechanical Engineering, University of Peradeniya, KY 20400, Sri Lanka\\ 
}

\IEEEauthorblockA{\IEEEauthorrefmark{3}Department of Mechanical Engineering, Texas Tech University, TX 79409, USA\\ 
}
}

\begin{document}

\maketitle
\begin{abstract}
This paper applies geometric PID control for asymptotic tracking of a desired trajectory by a hoop robot in the presence of disturbances and uncertainties. The hoop robot, consisting of a circular body rolling without slip along a one-dimensional surface, is a planar analog of a spherical robot. A variety of coupled mechanical system may be used to actuate the hoop robot. This paper specifically considers two different actuators, one a simple pendulum and the other an internal cart. The geometric PID controller requires the plant to be a mechanical system, and the hoop robot does not satisfy this condition. Therefore a geometric inner loop is presented that gives the hoop robot the required structure. This procedure is here referred to as feedback regularization. Feedback regularization---in contrast to feedback linearization---is coordinate independent, and hence reflects the fundamental system structure. Note also that the resulting mechanical system is nonlinear and underactuated. Subsequently, the geometric PID outer loop guarantees almost-semiglobal tracking with locally exponential convergence, and the integral action of the PID guarantees robustness to constant disturbances and parameter uncertainties, including constant inclination of the rolling surface. The complete tracking controller is the composition of the two coordinate-independent loops, and therefore is also coordinate independent.
\end{abstract}
\allowdisplaybreaks

\section{Introduction} \label{sect:Introduction}
Spherical robots combine a perfectly round body, rolling without slip on a planar surface, with an internal or external actuation mechanism. Spherical robots are a natural application for geometric control, because of the well-known difficulties in using a single set of coordinates to describe the dynamics of large rotations \cite{Mayhew}. The geometric approach  uses control formulations that remain valid in all coordinate systems, and so can be easily and consistently implemented in any convenient set \cite{Koditschek, BulloTracking, MaithripalaLieGroup, MaithripalaAutomatica}. Despite its apparent simplicity, the interactions between the spherical body, the constraint forces, and the dynamics of the actuation mechanism can lead to complex behavior and a substantial control challenge. This paper presents an approach to a feedback tracking controller for a class of spherical robots, that allows semi-global, locally exponential trajectory tracking in the presence of uncertainties, including a constant, unknown surface inclination.

Many existing results on spherical robots strictly concern open-loop path planning strategies \cite{Agrawal,Jurdjevic,Marigo,Mukherjee} or are valid only for perfectly horizontal surfaces \cite{Banavar2015,Cai,Karavaev}. To our knowledge, only one study \cite{Zhao} takes into account the inclination of the rolling surface. This controller combines feedback linearization with sliding mode control \cite{Zhao}. This controller of \cite{Zhao} is formulated in a single coordinate patch, and hence convergence is only guaranteed to be local. The controller of \cite{Zhao} also requires perfect knowledge of the inclination of the rolling surface. 

\begin{figure}[h!]
	\centering
		\begin{tabular}{cc}
		a) \includegraphics[width=0.15\textwidth]{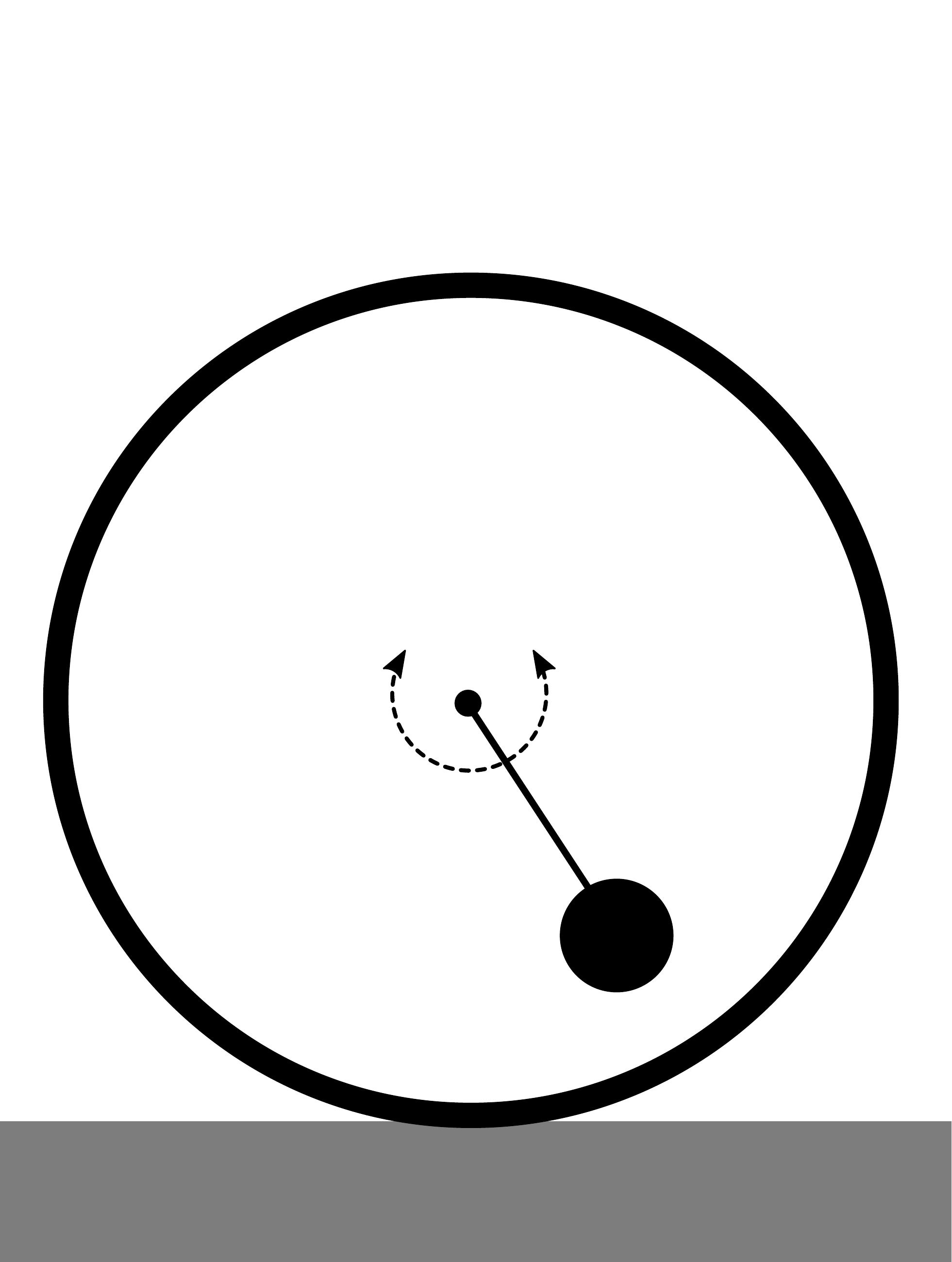}&
		b) \includegraphics[width=0.15\textwidth]{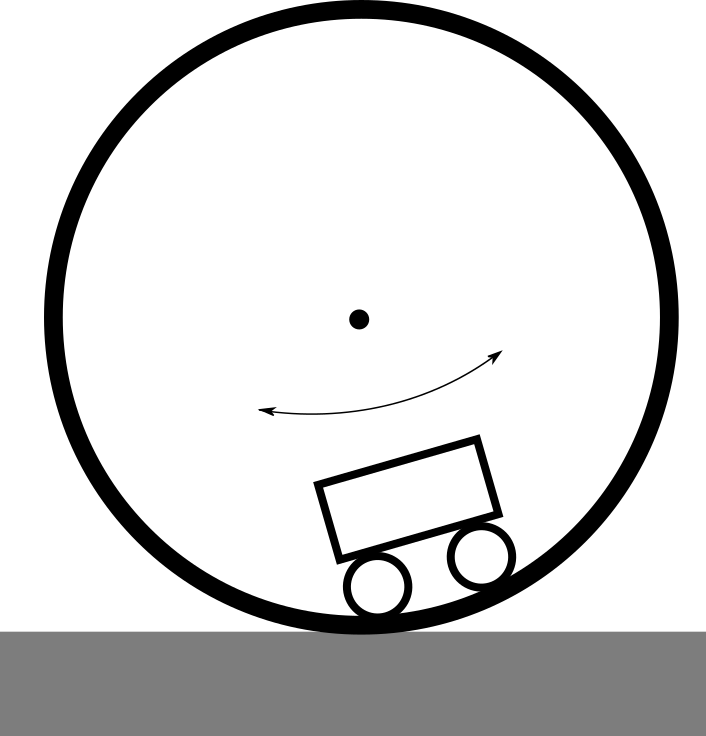}
	\end{tabular}
	\caption{Examples of hoop robots satisfying the constant-distance assumption.\label{Fig:BallbotYes}}
\end{figure}
\begin{figure}[h!]
	\centering
	\begin{tabular}{cc}
		a) \includegraphics[width=0.15\textwidth]{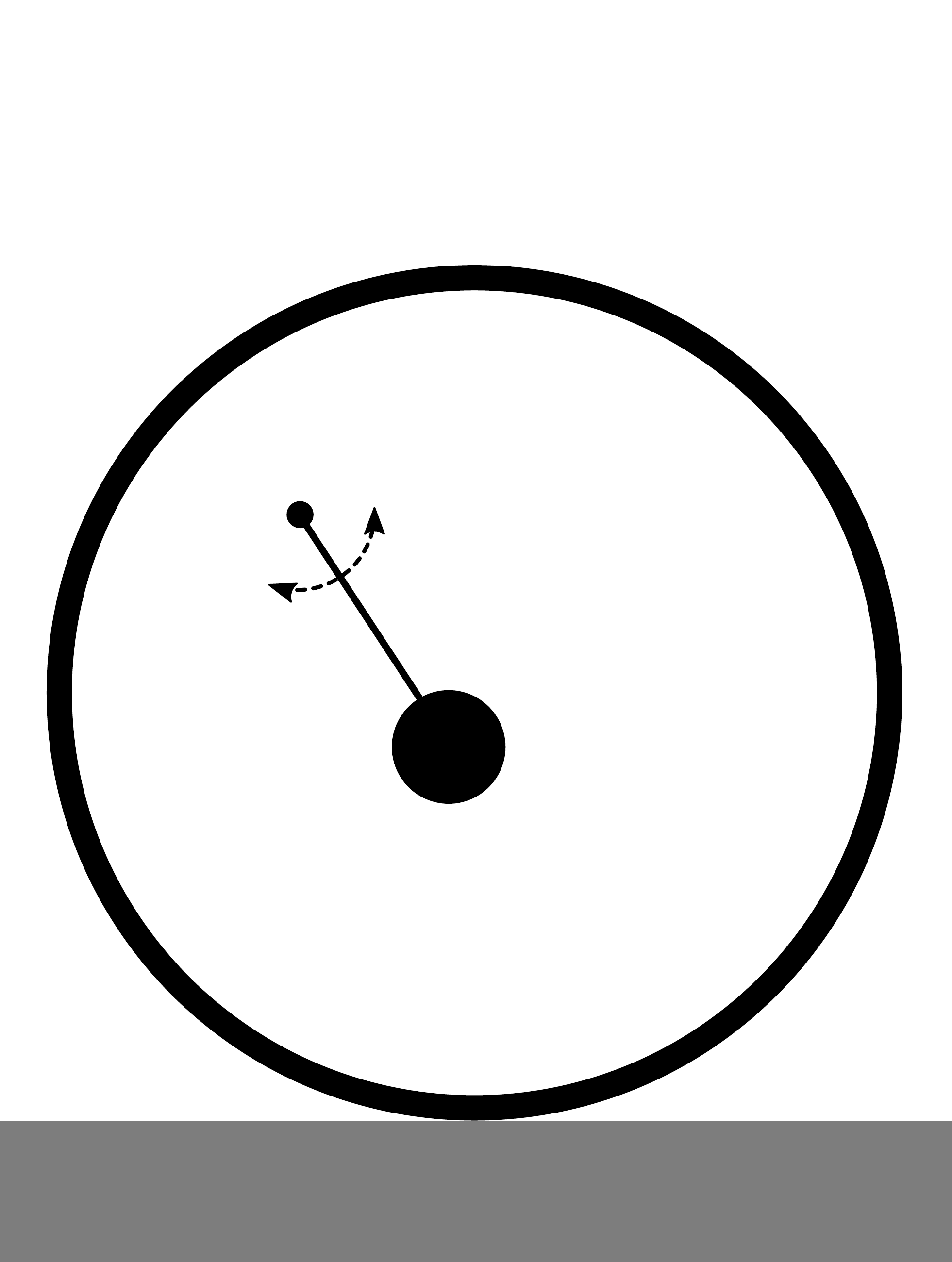}&
		b) \includegraphics[width=0.15\textwidth]{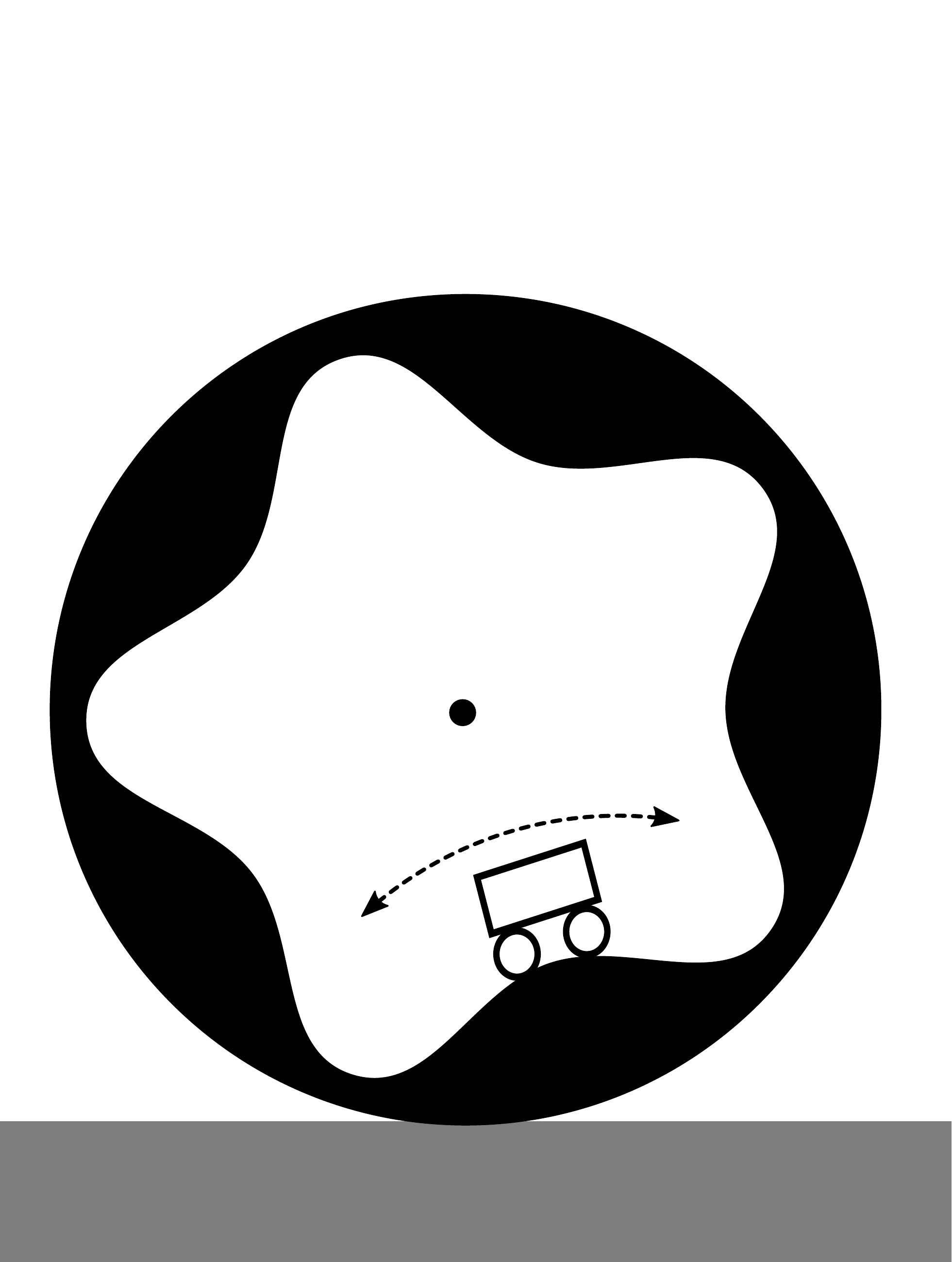}
	\end{tabular}
	\caption{Examples of hoop robots not satisfying the constant-distance assumption.\label{Fig:BallbotNo}}
\end{figure}
For simplicity of presentation, this approach is demonstrated on a class of \emph{hoop robots}---that is, on a planar version of a spherical robot. A hoop robot consists of a circular body, rolling without slip on a one-dimensional surface, and actuated by one of a variety of mechanisms. The class of hoop robots considered here requires that the center of mass of the moving actuation mechanism remain at constant distance from the geometric center of the hoop. This assumption is not fundamental; it is invoked to allow the actuator configuration space to be the circle, $\mathbb{S}$. Figure \ref{Fig:BallbotYes} shows examples of allowable mechanisms, including a simple pendulum driven around the hoop center, and a cart that follows a circular path around the hoop center. Figure \ref{Fig:BallbotNo} shows examples of mechanisms that are not considered in the subsequent analysis, including a simple pendulum driven around any point other than the hoop center, and a cart that follows a non-circular path. 

This analysis extends two familiar standard notions. The first is proportional + integral + derivative, or PID, control. The concept of proportional + derivative, or PD, control was first extended to the geometric setting in \cite{Koditschek}, and subsequently further developed in \cite{BulloTracking, MaithripalaLieGroup}. The central object in PD control is the tracking error, and the insight of the geometric extension is to give the tracking error dynamics the form of a \emph{mechanical system}. Such systems arise naturally in physics, and geometric PD control design follows from deriving control action from artificial ``potential energy'' terms that create ``energy minima'' where the tracking error is zero. While the ideas of geometric PD control can be applied to general Riemannian manifolds, the expressions are much more compact for systems evolving on Lie groups. More recently, geometric PD control has been augmented by integral action to obtain true geometric PID control \cite{MaithripalaAutomatica}. The key insight in this extension is that the natural form of differentiation of velocity for Riemannian manifolds is given by the Levi-Civita connection. Implementing an integral error term based on this insight,  \cite{MaithripalaAutomatica} obtains configuration tracking of a mechanical system on a Lie group. The integral action provides almost-global, locally exponential convergence of the tracking error to zero in the presence of bounded parametric uncertainty and bounded constant disturbance forces. However the result in \cite{MaithripalaAutomatica} only applies to  configuration tracking of fully-actuated systems. The broader class of systems that we consider here corresponds to two or more mechanical systems coupled through shared control forces. Considered independently, each subsystem is fully actuated, but since inputs are shared by two or more configuration variables, the combined system is underactuated. The subsystems may also be coupled through quadratic velocity terms arising from the Reimannian structure corresponding to the overall system kinetic energy. \emph{The first contribution of the present paper is to extend geometric PID control to output tracking of this class of underactuated mechanical systems}. 

Because the system is underactuated, it will not be possible to achieve configuration tracking of all configuration variables. Thus we designate one of the mechanical systems as the \emph{output system}, and collectively refer to the others as the \emph{actuator system}. We will be interested in ensuring that a output that depends only on the configuration of the output system tracks a desired output. We will assume that the output system is fully actuated and that the zero dynamics of the system has a stable but not necessarily  asymptotically stable equilibrium.
The \emph{output error system} is then the system describing the discrepancy between the output of the system and the desired output reference trajectory. The goal of the underactuated trajectory tracking problem is to achieve asymptotic convergence of the output error to zero while ensuring that the velocities of the actuator system remain bounded. 

While the class of coupled systems described above may seem overly narrow, it arises naturally when multi-body systems interact with each other. In that case the inputs are shared owing to the fact that the interaction forces and moments must be equal and opposite. 
Due to these interactions the expression of the dynamics of 
the individual subsystems may fail to correspond to simple mechanical systems. However we may use feedback control to make them look like interconnected simple mechanical systems. In some sense, this is a geometric interpretation of partial feedback linearization (PFL). PFL is a powerful technique of nonlinear control, in which state feedback is used to cancel all nonlinearities in the system input-output response \cite{Isidori}. In contrast rather than cancel terms we introduce terms that are quadratic in the velocities so that the individual systems will take the form of simple mechanical systems. Unlike the PFL procedure, our feedback terms (the terms we add) are independent of coordinate systems, and therefore can be used to provide the best possible global stability results. \emph{The second contribution of the present paper is to use coordinate-independent feedback to inject quadratic velocity terms that correspond to the Levi-Civita connection for the system kinetic energy and thereby provide each of the subsystems with the structure of a simple mechanical system}. Since the objective is not a linear system, but rather a simple mechanical system, we refer to this process as \emph{feedback regularization}. 
Under the assumptions of this paper, every system to which feedback regularization can be applied can subsequently be controlled to track a desired output trajectory using geometric PID. \textit{Therefore the third contribution of the current paper is to achieve general output trajectory tracking for the class of underactuated mechanical systems described above.}

Hoop robots provide an ideal test bed to demonstrate these techniques. The coupling between the hoop body and the actuation mechanism is through reaction forces and moments. If no other forces were present, Newton's laws applied to the error system augmented by coupled system would produce equations of motion suitable for the geometric PID tracking controller. However, the forces corresponding to the interaction between the subsystems destroy that structure. Thus we first use feedback regularization to recover the Riemannian structure of the subsystems, and subsequently apply geometric PID to the regularized system. 

In Section \ref{Secn:InterConnectedMechanicalSysms} we present our main results in the very general setting of interconnected mechanical systems on compact Riemannian manifolds where we extend the geometric PID controller of \cite{MaithripalaAutomatica} to simultaneously ensure almost-semi-global, locally exponential tracking of the output of a class of coupled underactuated mechanical systems. Appendix \ref{Secn:MechanicalSystems} contains a brief summary of helpful background material on the Riemannian geometry of mechanical systems.
%%%%%%%%%%%%%%
In Section \ref{sect:RollingHoop} we derive equations of motion for the hoop plus the actuation mechanism under the constant distance assumption. We see that equations that describe each subsystem do not have the structure of a mechanical system. In Section \ref{Secn:TrackingController} we show how to use feedback to give the coupled dynamic system consisting of the tracking error and the actuation mechanism the structure of an interconnected underactuated mechanical system on the circle, $\mathbb{S}$.  Section \ref{Secn:TrackingController} applies the general controller developed in Section \ref{Secn:InterConnectedMechanicalSysms} to ensure that the geometric center of the rolling hoop follows a desired trajectory, with boundedness of the velocities guaranteed for the actuation mechanism. Also in this section, the convergence properties of the controller are shown to be robust to bounded parametric uncertainties and constant unknown bounded disturbances, including those due to constant but unknown inclination.  We have extended these ideas to the full 3D sphere in \cite{RollingBallAutomatica}. Section \ref{sect:SimResults} presents simulations showing excellent performance for the complete system. 
%%%%%%%%%%%%%%%%%%%%%%%%%%%%%%

\section{Interconnected Mechanical Systems}\label{Secn:InterConnectedMechanicalSysms}
In the following we consider a class of interconnected mechanical systems on Riemannian manifolds. We refer the reader to \cite{Marsden,Bullo,Stocia} for further details on mechanical systems on Riemannian manifolds.

We will assume that each subsystem evolves on a configuration space $G_\nu$ and has an inertia tensor $\mathbb{I}_\nu$. Here $\nu$ is either $s$ or $a$, denoting the corresponding subsystem of the interconnected system. Denote by $\dot{g}_{\nu}=v_\nu$ and by $\nabla^\nu$ be the unique Levi-Civita connection corresponding to $\mathbb{I}_\nu$.

The class we consider is of the form
\begin{align}
\mathbb{I}_s\nabla^s_{v_s}{v_s} &=\tau_{s}(v_s,v_a)+\Delta_s+\tau_u,\label{eq:System_s}\\
\mathbb{I}_a\nabla^a_{v_a}{v_a} &=\tau^a_P(g_a)+\tau_{a}(v_s,v_a)+\Delta_a+B(g_a)\tau_u,\label{eq:System_a}
\end{align}
where $\Delta_s,\Delta_a$ represent input disturbances and unmodelled forces acting on the system, and $\tau_{s}$ and $\tau_{a}$ are quadratic velocity dependent interaction forces that satisfy $\tau_s(0,v_a(t))\equiv 0$, $\tau_a(0,v_a(t))\equiv 0$ for any $v_a(t)\in TG_a$. 
The system denoted by $s$ will be referred to as the \textit{ output system} and the system denoted by $a$ will be referred to as the \textit{actuation system}. The output system will be assumed to be fully actuated with respect to the input $\tau_u$.
Let $y: G_s \mapsto G_y$, %\label{eq:Output}
be a smooth onto function for some smooth Lie-group $G_y$.  We will also assume that $y$ is relative degree two with respect to $\tau_u$. The output that we are interested in will be $(g_y(t),v_s(t))\in G_y\times TG_s$.
The control problem that we solve in this section is that of ensuring the almost-semi-global and local exponential convergence of $(g_y(t),v_s(t))$ to $(e_y,0)$ where $e_y$ is the identity element of $G_y$. 

From (\ref{eq:System_s}) it follows that the output zeroing control $\bar{\tau}_u$ must necessarily satisfy
\begin{align}
\bar{\tau}_u&=-\Delta_s.\label{eq:Eqms}
\end{align}
The corresponding zero dynamics of the system is given by
\begin{align}
\mathbb{I}_a\nabla^a_{v_a}{v_a} &=\tau^a_P(g_a)+\Delta_a-B(g_a)\Delta_s.\label{eq:GeneralZeroDy}
\end{align}

We will assume that for any given bounded constant disturbances the zerodynamics (\ref{eq:GeneralZeroDy}) have a stable relative equilibrium.
We will state this specifically in the following assumption:
\begin{assumption}\label{as:AssumptionInterconnected}
For any given bounded constant $\Delta_s,\Delta_a$ there exists a smooth positive semi-definite potential function $V_a: G_a\mapsto \mathbb{R}$ such that 
\[
dV_a=-(\tau^a_P(g_a)+\Delta_a-B(g_a)\Delta_s).
\]
\end{assumption} 
What this assumption implies is that ensuring $\lim_{t\to \infty}(g_y(t),v_s(t))=(e_y,0)$ with $\lim_{t\to \infty}\tau_u(t)=-\Delta_s$ exponentially guarantees that the velocity of the actuator system $v_a(t)$ remains bounded.

%%%%%%%%%%%%%%%%%%%%%%%%%%%%%%%%%%%%

\subsection{Nonlinear PID Control for Interconnected Mechanical Systems}\label{Secn:NLPID}
The control problem that we solve in this section is that of ensuring the almost-semi-global and local exponential convergence of the smooth output $(y(t),v_s(t))$ to $(e_y,0)$ where $e_y$ is the identity on $G_y$.
Let $V_{y} : G_y \mapsto \mathbb{R}$ be a polar Morse function on $G_y$ with a unique minimum at $y_0$. Let $\eta_s\in T_{g_s}G_s$ be the pullback of the gradient of $V_{y}$ that is defined by $\langle\langle\eta_s,v_s\rangle\rangle=dV_s(v_s)$ where $V_s(g_s)\triangleq V_y(y(g_s))$.

Consider the nonlinear potential shaping plus PID controller
\begin{align}
&\mathbb{I}_s\nabla^s_{v_s} {v_I} =\mathbb{I}_s{\eta_e},\label{eq:GeneralI}\\
&\tau_u=-\mathbb{I}_s(k_{p}{\eta_s}+k_{d}{v_s}+k_{I}{v_I}),\label{eq:GeneralPID}
\end{align}
where $v_I,\eta_e\in T_{g_s}G_s$ with the gains $k_p,k_I,k_d$ chosen such that
\begin{align}
0<k_I<\frac{k_d^3(1-\delta^2)}{\mu},\label{eq:kICond}\\
 k_p>\max \left\{k_1,k_2,2\kappa k_d^2\right\},\label{eq:kpCond}
 \end{align}
 where,
 \begin{align*}
 k_1&=\frac{k_I}{2k_d}\left(\sqrt{1+\frac{16r\kappa^2k_d^2}{k_I}}-1\right),\\
 k_2&=\frac{ r k_I^2}{2k_d^4}\left({1 +\sqrt{1+\frac{4k_d^3 (k_I^2+4\kappa k_d^3(1+\kappa k_d^3))}{r k_I^3}}}\right).
\end{align*}

In the Appendix we prove the following theorem.
\begin{theorem}\label{theom:Theorem1} 
Consider any compact subset $\mathcal{X}_s\subset G_y\times \mathcal{G}_s\times \mathcal{G}_s$ that contains $(0,0,0)$.
Assume that the conditions of Assumption-\ref{as:AssumptionInterconnected} hold. Then if the gains of the nonlinear PID controller (\ref{eq:GeneralI})--(\ref{eq:GeneralPID}) are chosen to satisfy (\ref{eq:kICond}) and (\ref{eq:kpCond}) the followings hold for almost all initial conditions in $\mathcal{X}_s$:
\begin{enumerate}
\item $(y(t),v_s(t),v_I(t))$ converges to an arbitrarily small neighborhood of $(y_0,0,0)$ while,
\item if the uncertainties and the disturbances are constant then 
$\lim_{t\to \infty}(y(t),v_s(t))=(y_0,0)$ locally exponentially.
\end{enumerate} 
\end{theorem}

%%%%%%%%%%%%%%%%%%%%%%%%%%%%%  

\section{Planar Hoop Dynamics on an Inclined Ramp} \label{sect:RollingHoop}
\begin{figure}[h!]
	\centering
	\begin{tabular}{c}
		\includegraphics[width=0.30\textwidth]{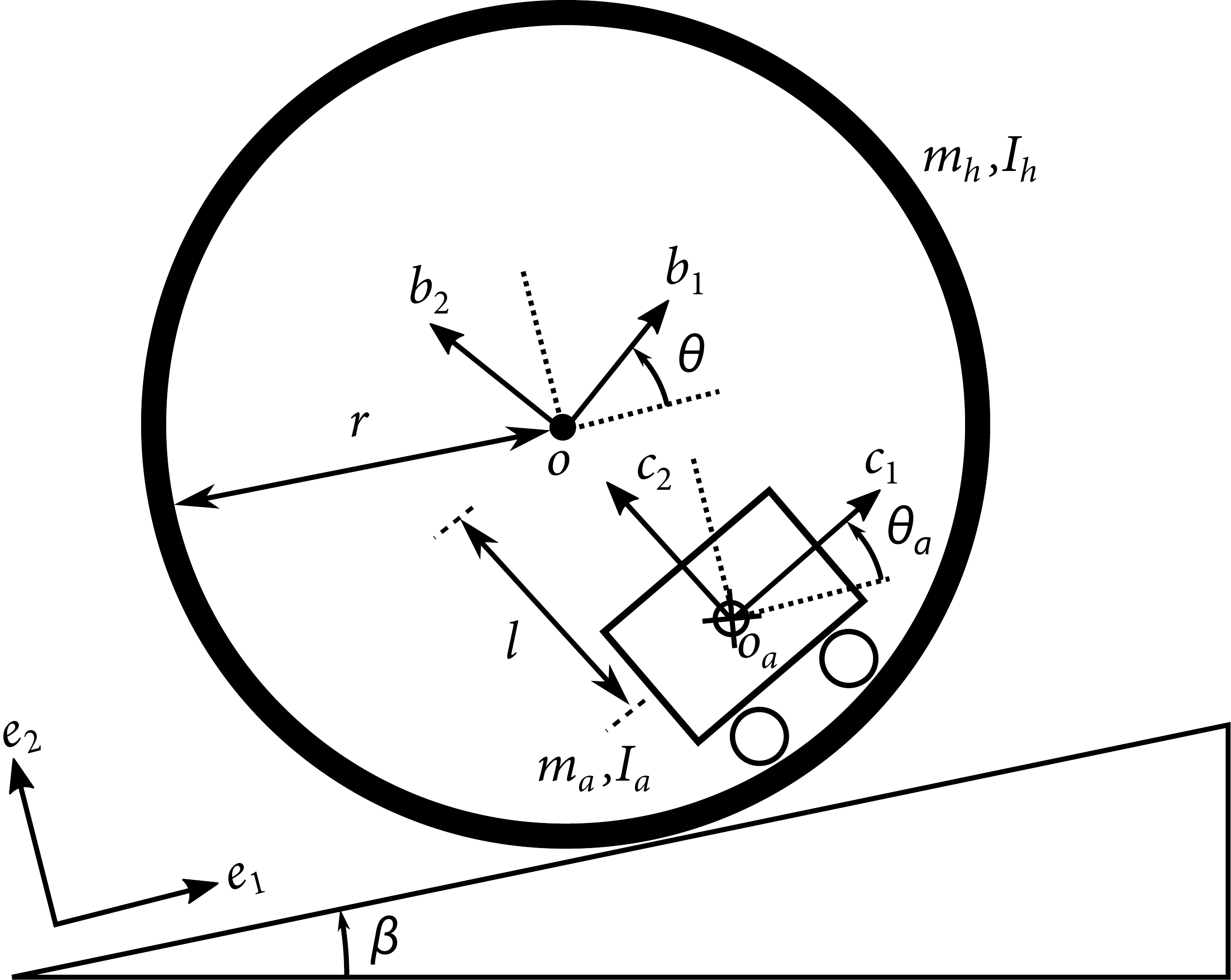}
	\end{tabular}
	\caption{Reference frames for geometric hoop robot analysis. The frames are well-defined for any actuation mechanism satisfying the constant-distance assumption.\label{Fig:HoopRefFrames}}
\end{figure}
Figure \ref{Fig:HoopRefFrames} shows a hoop of radius $r$, rolling without slip on an inclined, one-dimensional surface of constant inclination $\beta$ with respect to the horizontal plane. Let $\mathbf{e}=[\mathbf{e}_1\:\:\:\mathbf{e}_2]$ be an earth-fixed inertial frame with second axis pointing along the outward normal of the surface. Let $\mathbf{b}=[\mathbf{b}_1\:\:\:\mathbf{b}_2]$ be a reference frame fixed on the rolling hoop, with origin coinciding with the geometric center of the hoop. Let $\theta$ be the angle of rotation of the frame $\mathbf{b}$ with respect to the frame $\mathbf{e}$, and let $\omega=\dot{\theta}$. Let $m_h$ be the mass of the hoop and $\mathbb{I}_h$ be the inertia of the hoop. We assume that a mechanism of mass $m_a$ and moment of inertia $\mathbb{I}_a$ actuates the hoop. Let $\mathbf{c}=[\mathbf{c}_1\:\:\:\mathbf{c}_2]$ be a frame fixed to the actuation mechanism, with origin $o_a$ coinciding with the actuation center of mass and axis $\mathbf{c}_2$ pointing towards the center of the hoop. Let $\theta_a$ be the rotation angle of $\mathbf{c}$ with respect to $\mathbf{e}$ and $\omega_a=\dot{\theta}_a$.

Let $l$ be the distance from the geometric center of the hoop to the center of mass of the actuation mechanism. We restrict our attention to cases where $l$ remains constant. Under this assumption, the actuation mechanism evolves on the circle $\mathbb{S}$, and is completely characterized by configuration variable $\theta_a$. While this assumption is not strictly necessary to the rest of the analysis, it allows a compact characterization of the result.

We assume that the external forces and external moments acting on the hoop are due to the effect of gravity, reactions that arise as a consequence of the interaction with the actuation mechanism and the no-slip rolling constraint. The gravity force acting on the hoop is $f_g=-m_hg\,e_g$ where $e_g$ is the unit vector in the vertical direction, which can be written as $\sin\beta\,  \mathbf{e}_1 + \cos\beta\,  \mathbf{e}_2$ with respect to $\mathbf{e}$. Let $f_{\lambda}=[f_{{\lambda}_1} \:\:\: f_{{\lambda}_2}]^T \in \mathbb{R}^2$ be the $\mathbf{e}$-frame representation of the force that ensures the no-slip rolling constraints. The resultant moment acting on the hoop due to these constraint forces is $\tau_{\lambda}=rf_{{\lambda}_1}$. Let $f_c=[f_{c_1} \:\:\: f_{c_2}]^T \in \mathbb{R}^2$ be the $\mathbf{e}$-frame representation of the force acting on the hoop due to the interaction with the actuation mechanism and let $\tau_{f_c}$ be the resultant moment acting on the hoop. All moments are assumed to be taken with respect to the geometric center of the hoop. By Newton's third law, $-f_c$ and $-\tau_{f_c}$ are the reaction forces and moments acting on the actuation mechanism. For cart-type actuation mechanism the interaction between the hoop and the cart occurs only through no-slip constraints at the cart wheels, and hence $\tau_c\equiv 0$. For pendulum-type actuation, $\tau_{f_c}\equiv 0$.

For the hoop, Euler's rigid body equations and the hoop/ground no-slip constraint give
\begin{align*}
\mathbb{I}_h\dot{\omega}&=\tau_c+\tau_{f_c}+rf_{\lambda_1}.\\
f_{\lambda_1}&=-m_hg\sin\beta-f_{c_1}-m_hr\dot{\omega},
\end{align*}
\noindent and hence 
\begin{align*}
(\mathbb{I}_h+m_hr^2)\dot{\omega}=rm_hg\sin\beta-rf_{c_1}+\tau_c+\tau_{f_c}.
\end{align*}

For the actuator, Euler's rigid body equations give 
\begin{align*}
\mathbb{I}_a\dot{\omega}_a&=l\sin\theta_af_{c_2}+l\cos\theta_af_{c_1}
-(\tau_c+\tau_{f_c}),\\
f_{c_1}&=
-m_ag\sin\beta+m_a(r\dot{\omega}-l\dot{\omega}_a\cos\theta_a+l\omega_a^2\sin\theta_a)\\
f_{c_2}&=-m_ag\cos\beta-m_al(\dot{\omega}_a\sin\theta_a+\omega_a^2\cos\theta_a).
\end{align*}

Defining
{\tiny
\begin{align*}
M&\triangleq m_h+m_a , \\
\mathbb{I}(\theta_a)&\triangleq\left(\mathbb{I}_h+Mr^2-\frac{m_a^2r^2l^2}{\mathbb{I}_a+m_al^2}\cos^2{\theta_a}\right),\\
\tau^\omega_g&\triangleq rMg\sin{\beta}-\frac{m_a^2rl^2g}{\mathbb{I}_a+m_al^2}\cos{\theta_a}\sin{(\theta_a+\beta)},\\
\tau^{\omega_a}_g&\triangleq \frac{m_arl\cos{\theta_a}}{\mathbb{I}_a+m_al^2}\tau^\omega_g-\mathbb{I}(\theta_a)\left(\frac{m_agl\sin({\theta_a+\beta})}{\mathbb{I}_a+m_al^2}\right) , \\
B(\theta_a)&\triangleq \left(\frac{m_arl\cos{\theta_a}}{\mathbb{I}_a+m_al^2}-\frac{\mathbb{I}(\theta_a)}{(\mathbb{I}_a+m_al^2-m_arl\cos{\theta_a})}\right),
\end{align*}} 
\noindent yields the following complete hoop robot equations of motion on the state space $\mathbb{S}\times\mathbb{R}\times\mathbb{R}\times\mathbb{S}\times\mathbb{R}$:
{\small
	\begin{align} 
	&\dot{\theta}=\omega,\\
	&\dot{o}=-r\omega,\\
	&\mathbb{I}(\theta_a)\dot{\omega}=-m_arl\sin{\theta_a}\omega_{a}^2+\tau^\omega_g+\tau_u,\\
	&\dot{\theta}_a={\omega}_a,\\
&\mathbb{I}(\theta_a)\dot{\omega}_a=
-\frac{m^2_ar^2l^2\sin{\theta_a}\cos{\theta_a}}{\mathbb{I}_a+m_al^2}\omega_{a}^2+\tau^{\omega_a}_g+B(\theta_a)\tau_u .  \label{eq:CartDy}
\end{align}}

The single control input, which appears in both the $\omega$ and $\omega_a$ equations, is defined as
\begin{align*}
\tau_u\triangleq \frac{(\mathbb{I}_a+m_al^2-m_arl\cos{\theta_a})}{\mathbb{I}_a+m_al^2}\left(\tau_c+\tau_{f_c}\right) .
\end{align*}

%%%%%%%%%%%%%%%%%%%%%%%%%%%%%%%%%%%%%%%%%%%%%%
\section{Position Tracking for the Hoop}\label{Secn:TrackingController}

The control task that we consider is to ensure that the output satisfies $\lim_{t \to \infty} o(t)=o_{\mathrm{ref}}(t)$ where $o_{\mathrm{ref}}(t)\in \mathbb{R}$ is a twice differentiable reference  and $\omega_{\mathrm{ref}}=-\dot{o}_{\mathrm{ref}}/r$. Let $o_e\triangleq (o-o_{\mathrm{ref}})$, $\omega_e\triangleq (\omega-\omega_{\mathrm{ref}})$. Note that in the special case of stabilizing the hoop at a point the references are constant: $o_{\mathrm{ref}}(t)\equiv \mathrm{const}$ and $\omega_{\mathrm{ref}}\equiv 0$.

Differentiating these quantities one sees that the error dynamics of the system take the explicit form
\begin{align}
\dot{o}_e&=-r\omega_e,\label{eq:Erroro}\\
\mathbb{I}({\theta_a})\dot{\omega}_e&=-m_arl\sin{\theta_a}\omega_{a}^2+\tau_g^{\omega}-\tau_{\mathrm{ref}}+\tau_u,\label{eq:Erroromega}
\end{align}
where $\mathbb{I}(\theta_a)\triangleq \mathbb{I}_h+Mr^2-\frac{m_a^2r^2l^2}{\mathbb{I}_a+m_al^2}\,\cos^2\theta_a$, and
$\tau_{\mathrm{ref}}\triangleq\mathbb{I}(\theta_a)\dot{\omega}_r$.
We notice that the error dynamics evolve on the tangent bundle to the circle, $T\mathbb{S}$, with output $y=o_e$ evolving on the Lie-group $\mathbb{R}$.
The natural notion of differentiation on the tangent space of a Riemannian manifold is the Levi-Civita connection, $\nabla$. As discussed in Appendix \ref{Secn:MechSysOnS}, the unique Levi-Civita connection on $\mathbb{S}$ corresponding to the kinetic energy induced by the inertia tensor $\mathbb{I}$ is explicitly given by
\begin{align*}
\mathbb{I}({\theta_a})\nabla_{\zeta}\eta=\mathbb{I}({\theta_a})d\eta(\zeta)+\frac{m_a^2r^2l^2\sin{(2\theta_a)}}{2(\mathbb{I}_a+m_al^2)}\zeta\eta,
\end{align*}
for any $\zeta,\eta\in T_{g_a}\mathbb{S}$. Setting $\zeta$ to $\omega_a$ and $\eta$ to $\omega_e$ gives
\begin{align}
\mathbb{I}({\theta_a})\nabla_{\omega_a}\omega_e = \mathbb{I}({\theta_a})\dot{\omega}_e + \frac{m_a^2r^2l^2\sin{(2\theta_a)}}{2(\mathbb{I}_a+m_al^2)}\omega_a \omega_e.\label{eq:LeviCivitaS}
\end{align}

As discussed in Appendix \ref{Secn:MechanicalSystems}, using the Levi-Civita connection, the equations of motion of a mechanical system are
\begin{align}
\mathbb{I}\nabla_{\dot{g}}\dot{g}&=\gamma_g . \label{eq:MechanicalSystemm}
\end{align} 
where $\gamma_g$ is the generalized force acting on the system. Comparing with (\ref{eq:LeviCivitaS}) and (\ref{eq:MechanicalSystemm}) it is clear that , because of the absence of the term
\begin{align*}
\frac{m_a^2r^2l^2\sin{(2\theta_a)}}{2(\mathbb{I}_a+m_al^2)}\omega_a \omega_e
\end{align*}
(\ref{eq:Erroromega}) does not represent a simple mechanical system.
This prevents the straightforward use of the nonlinear PID controller proposed by the authors in \cite{MaithripalaAutomatica}.

However we notice that if we choose the regularizing plus potential shaping controls
\begin{multline}
\tau_u=-\frac{m_a^2r^2l^2\,\sin{2\theta_a}}{2(\mathbb{I}_a+m_al^2)}\,\omega_a\omega_e-m_arl\sin{\theta_a}\omega_{a}^2\\+\frac{m_a^2rl^2g\sin{(2\theta_a)}}{2(\mathbb{I}_a+m_al^2)}+\tilde{\tau}_u,\label{eq:TransformingControl}
\end{multline}
\noindent then the error dynamics of the system (\ref{eq:Erroromega}) can be re-written as,
\begin{align}
\mathbb{I}({\theta_a})\nabla_{\omega_a}\omega_e&=\Delta_h+\tilde{\tau}_u.\label{eq:IntrinsicOmegaeDy}
\end{align}
The third term in the above equation (\ref{eq:TransformingControl}) shapes the potential energy of the error dynamics.
Here $\Delta_h$ represents the effects due to the ignorance of the inclination of the rolling surface and the omission of the term $\tau_{\mathrm{ref}}$ in the controller.
In similar fashion we find that the actuation system dynamics can also be expressed as
{\small
\begin{align}
\mathbb{I}\nabla_{\omega_a}\omega_a&=\tilde{\tau}^{\omega_a}_g+\Delta_a+\tau_{a}(\omega_e,\omega_a)+B(\theta_a)\tilde{\tau}_u.\label{Eq:InnerCart}
\end{align}}
\noindent where $\Delta_a$ represents modeling errors and disturbances. Here 
\begin{align*}
\tilde{\tau}^{\omega_a}_g&\triangleq {\tau}^{\omega_a}_g +B(\theta_a)\left(\frac{m_a^2l^2g\sin{(2\theta_a)}}{2(\mathbb{I}_a+m_al^2)}\right),\\
\tau_{a}(\omega_e,\omega_a)&\triangleq -B(\theta_a)\left(\frac{m_a^2r^2l^2\,\sin{2\theta_a}}{2(\mathbb{I}_a+m_al^2)}\,\omega_a\omega_e\right).
\end{align*}
Notice that $\tau_a(\zeta,\eta)$ is bilinear in the two velocity arguments $\zeta,\eta\in\mathbb{R}$. 
 When $\omega_e\equiv 0$ it can be shown that there exists an equilibrium for the actuation mechanism (\ref{Eq:InnerCart}) for any surface of inclination $\beta \in (-\pi/2,\pi/2)$ if the system parameters satisfy, $\frac{m_al}{Mr}\geq \sin\beta$. Without loss of generality we assume that the actuation mechanism is chosen such that this condition is satisfied for a certain operating region of $\beta$.

Note that the combination of the error dynamics of the system (\ref{eq:IntrinsicOmegaeDy}) and the actuation mechanism dynamics (\ref{Eq:InnerCart}) takes the form of the general interconnected under actuated mechanical system (\ref{eq:System_s})--(\ref{eq:System_a}) presented in Section \ref{Secn:InterConnectedMechanicalSysms}. The system evolves on $(\mathbb{S}\times \mathbb{R})\times (\mathbb{S}\times\mathbb{R})$ with the output $y=o_e$ evolving on the Lie-group $\mathbb{R}$.
For this interconnected mechanical system we propose the nonlinear PID controller,
\begin{align}
&\mathbb{I}({\theta_a})\nabla_{\omega_a}o_I=\mathbb{I}({\theta_a})\eta_e,\label{eq:IntegratorError}\\
&\tilde{\tau}_u=-\mathbb{I}({\theta_a})(k_p\eta_e+k_d\omega_e+k_Io_I),\label{eq:PIDcontroller}
\end{align}
where, $\eta_e=-o_e$.

Theorem \ref{theom:Theorem1} yields the following corollary:
\begin{corollary}
Assume that the parameter uncertainty, unmodelled disturbances, and the velocity references are bounded and constant.  There exists sufficiently large gains of the nonlinear PID controller (\ref{eq:IntegratorError})--(\ref{eq:PIDcontroller}) that satisfy (\ref{eq:kICond})--(\ref{eq:kpCond}) such that 
$\lim_{t\to \infty} (o_e(t),\omega_e(t))=0$ semi-globally and exponentially while ensuring that $\omega_a(t)$ remains bounded for any $(\theta_a(0),\omega_e(0),\omega_I(0),\omega_a(0))\in \mathcal{X}$, where $\mathcal{X}\subset \mathbb{S}\times\mathbb{R}\times\mathbb{R}\times \mathbb{R}$ is compact.  Here  $1/\mu<\kappa<2/\mu$ where $\mu=1+\frac{m_a^2r^2l^2k_\mathcal{X}}{2\sqrt{(\mathbb{I}_h+Mr^2)((\mathbb{I}_h+Mr^2)(\mathbb{I}_a+m_al^2)-m_a^2r^2l^2)}}$.
\end{corollary}

%%%%%%%%%%%%%%%%%%%%%%%%%%%%

\section{Simulation Results} \label{sect:SimResults}
In here we present simulation results that demonstrate the effectiveness of the nonlinear PID tracking controller for a rolling hoop on an inclined plane. The actuation mechanism we consider is of the type of a cart or a pendulum. In order to demonstrate the robustness of the controller we choose the system parameters in the simulations to be $50\%$ different from the nominal parameters used in the controller. 

 The nominal parameters for the hoop were chosen to be $m_h=1.00\, \si{\kg}$ and $\mathbb{I}_h=0.021\, \si{\kg.\m^2}$ while the outer radius of the hoop was chosen to be $r=0.18\, \si{\m}$. These parameters were chosen to correspond to a hoop made of plastic (density $850\, \si{\kg.\m^{-3}}$) with thickness of $3\, 
 \si{\mm}$. The nominal parameters for the actuation appendage was chosen to be $m_a =3.28\, \si{\kg}$ and $\mathbb{I}_a=0.035\, \si{\kg.\m^2}$. The distance from the geometric center of the hoop to the center of mass of the actuation mechanism was chosen to be $l=0.14\, \si{\m}$. For these parameters we find that the maximum inclination of the rolling plane for which an equilibrium exists for the actuation mechanism is $\beta_{\mathrm{max}}=36^\circ$.
Thus in the simulations the hoop is assumed to roll on a inclined plane of $20^\circ < \beta_{\mathrm{max}}$. We stress that the implementation of the controller (\ref{eq:IntegratorError})--(\ref{eq:PIDcontroller}) does not require the knowledge of the angle of inclination of the rolling surface.

\begin{figure}[h!]
	\centering
	\begin{tabular}{cc}
		\includegraphics[width=0.235\textwidth]{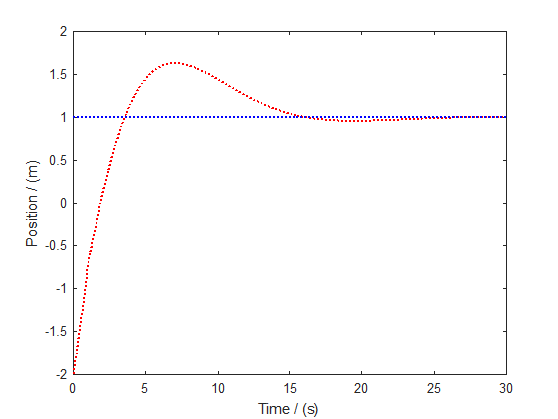}& \includegraphics[width=0.235\textwidth]{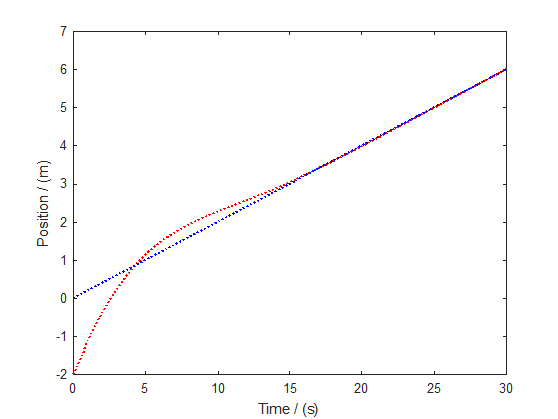}\\
		(a) Fixed point & (b) Linear velocity profile\\
		\includegraphics[width=0.235\textwidth]{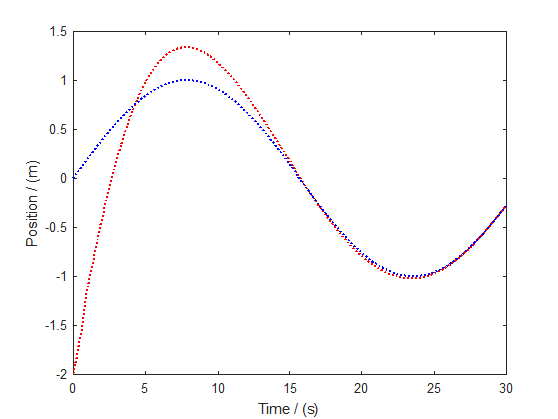}\\
		(c) Sinusoidal velocity profile
	\end{tabular}
	\caption{The position followed by the geometric center of the hoop for the PID controller (\ref{eq:IntegratorError})--(\ref{eq:PIDcontroller}) in the presence of parameter uncertainties as large as $50\%$. The blue curve shows the reference position while the red curve shows the actual position of the center of the hoop.\label{Fig:RollingHoopPosition}}
\end{figure}

\begin{figure}[h!]
	\centering
	\begin{tabular}{cc}
		\includegraphics[width=0.235\textwidth]{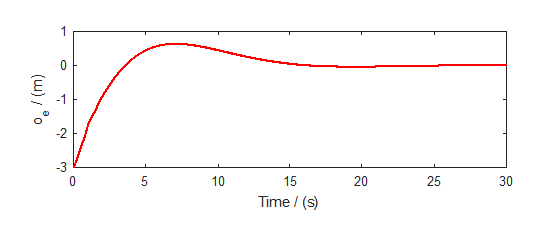}& \includegraphics[width=0.235\textwidth]{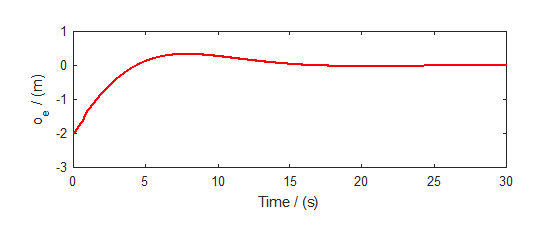}\\
		(a) Fixed point & (b) Linear velocity profile\\
		\includegraphics[width=0.235\textwidth]{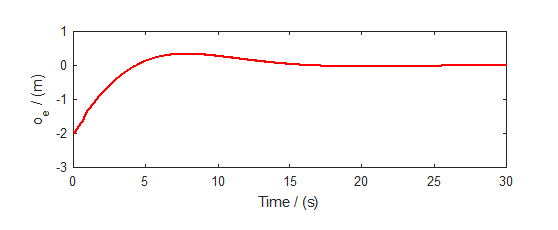}\\
		(c) Sinusoidal velocity profile
	\end{tabular}
	\caption{The position error, $o_e(t)$, for the PID controller (\ref{eq:IntegratorError})-(\ref{eq:PIDcontroller}) in the presence of parameter uncertainties as large as $50\%$.\label{Fig:RollingHoopPositionError}}
\end{figure}

\begin{figure}[h!]
	\centering
	\begin{tabular}{cc}
		\includegraphics[width=0.235\textwidth]{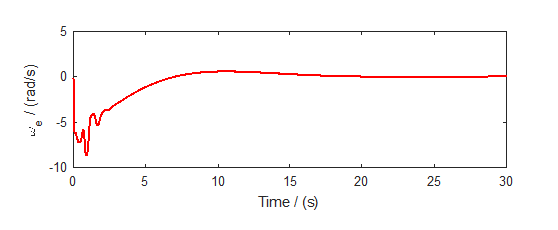}& \includegraphics[width=0.235\textwidth]{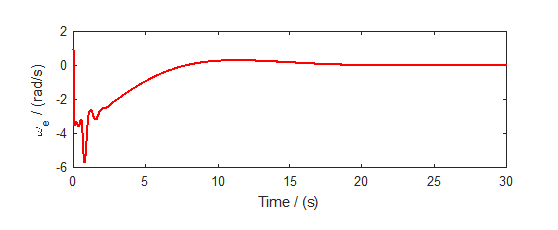}\\
		(a) Fixed point & (b) Linear velocity profile\\
		\includegraphics[width=0.235\textwidth]{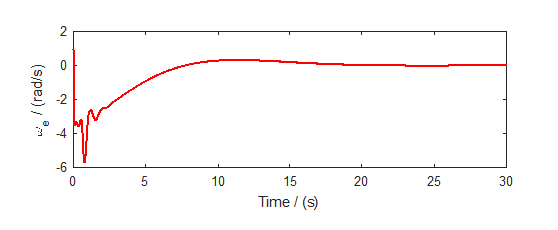}\\
	(c) Sinusoidal velocity profile
	\end{tabular}
	\caption{The spatial angular velocity error, $\omega_e(t)$, for the PID controller (\ref{eq:IntegratorError})--(\ref{eq:PIDcontroller}) in the presence of parameter uncertainties as large as $50\%$.\label{Fig:RollingHoopOmegaeError}}
\end{figure}

\begin{figure}[h!]
	\centering
	\begin{tabular}{cc}
		\includegraphics[width=0.235\textwidth]{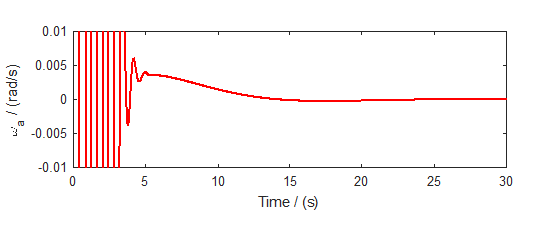}& \includegraphics[width=0.235\textwidth]{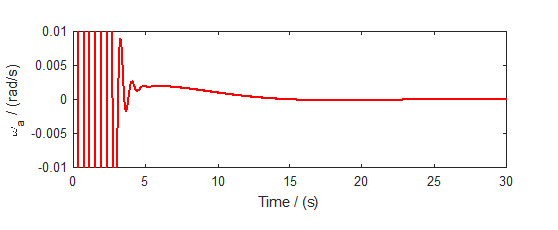}\\
		(a) Fixed point & (b) Linear velocity profile\\
		\includegraphics[width=0.235\textwidth]{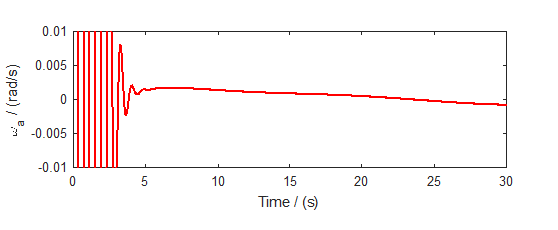}\\
		(c) Sinusoidal velocity profile
	\end{tabular}
	\caption{The spatial angular velocities of the  inner cart, $\omega_{a}(t)$, for the PID controller (\ref{eq:IntegratorError})--(\ref{eq:PIDcontroller}) in the presence of parameter uncertainties as large as $50\%$.\label{Fig:RollingHoopomegai}}
\end{figure}

We present results for: a.) stabilizing the hoop at a point, b:) tracking a linear velocity, and c.) tracking a sinusoidal velocity. The simulation results are presented in figure (\ref{Fig:RollingHoopPosition})--(\ref{Fig:RollingHoopomegai}). Figure-\ref{Fig:RollingHoopPosition} shows the reference position and actual position of the geometric center of the hoop, Figure-\ref{Fig:RollingHoopPositionError} and figure-\ref{Fig:RollingHoopOmegaeError} shows the position error of the geometric center of the hoop and the spatial angular velocity error of the hoop respectively while in figure-\ref{Fig:RollingHoopomegai} we illustrate that the spatial angular velocity of the actuation mechanisms remain bounded. Specifically, as expected, the actuator velocities tend to zero in the case of the constant set point and the constant reference velocity while it does not for the sinusoidal velocity profile.
In all simulations the initial position of the hoop was assumed to be at $o(0)=[-2\:\:\:r]^T\,\si{\m}$ and the initial spatial angular velocity for the hoop and inner cart were chosen as $\omega(0)=-0.1\,\si{\radian.\s^{-1}}$ and $\omega_a(0)=0.1\,\si{\radian.\s^{-1}}$ respectively. In all simulations the controller gains were chosen to be $k_p=16, k_d=7, k_I=4, k_c=0.1$.

%%%%%%%%%%%%%%%%%%%%%%%%%%%%% 

\section{Conclusions} \label{sect:Conclusions}
In this paper we present a control strategy for semi-almost-global output tracking for a class of interconnected under actuated mechanical systems. The control strategy involves two steps. In the first step, feedback control is used to give each of the subsystems the structure of a simple mechanical system. We call this feedback regularization. Next we use PID control to ensure that the output of one of the systems tracks desired trajectory while ensuring the stability of the other system. These results are the applied to tracking the geometric center of a hoop rolling without slip on an inclined plane of unknown inclination. The problem is a one dimensional approximation of a spherical robot. The controller is shown to ensure semi-almost-global exponential tracking in the presence of bounded parametric uncertainties and bounded constant disturbances. 
%%%%%%%%%%%%%%%%%%%%%%%%%%%%%y  

\begin{appendix}

%%%%%%%%%%%%%%%%%%%%%%%%%%%%%%%%%%%%%%
\subsection{Mechanical Systems on Riemannian Manifolds}\label{Secn:MechanicalSystems}
A mechanical system is defined by a configuration space, the kinetic energy, the generalized forces, and holonomic and non-holonomic constraints \cite{Marsden,Bullo}. 
The configuration space $G$ has the structure of a (possibly singular) Riemannian manifold. In most physical systems the $G$ is a smooth Lie-group. Denote by $T_gG$ the tangent space to $G$ at $g\in G$. The collection of all such tangent spaces to $G$ is referred to as the tangent bundle, $TG$ . The generalized velocity of the system at $g$, denoted by $v_g=\dot{g}$, is an element of the tangent space $T_gG$ while the generalized force acting on the system at $g$, denoted by $\gamma_g$, is an elements of the cotangent space $T_g^*G$ where $T_g^*G$ is the space dual to $T_gG$. It is customary to denote by $\langle\cdot,\cdot\rangle :
T_g^*G \times T_gG \mapsto \mathbb{R}$ the action of a covector $\gamma_g\in T_g^*G$ on a vector $v_g\in T_gG$ explicitly given by $\langle \gamma_g,v_g\rangle$. 

A singular Riemannian metric on $G$ assigns in a smooth fashion a degenerate inner product, 
$\langle\langle \cdot,\cdot \rangle \rangle : T_g G\times T_g G \mapsto \mathbb{R}$, on each of the tangent spaces $T_g G$ and gives $G$ the structure of a singular Riemannian manifold.  This metric is chosen such that the kinetic energy of the mechanical system is given by $\mathrm{KE}=\langle\langle \dot{g},\dot{g} \rangle \rangle/2$.
The singular Riemannian metric allows one to define a map $\mathbb{I}_g\::\:T_gG\mapsto T^*_gG$ by the relationship $\langle\mathbb{I}_gv_g,u_g\rangle\triangleq \langle\langle v_g,u_g\rangle\rangle$ for all $v_g,u_g\in T_gG$. The smooth tensor field $\mathbb{I}\::\:TG\mapsto T^*G$ that is point wise defined  above is usually referred to as the \textit{inertia tensor}.  For  Riemannian metrics the above defined map $\mathbb{I}_g$ is an isomorphism and thus in this case one can uniquely identify a vector with a given covector in an intrinsic fashion.

For any vectorfields $X,Y,Z\in TG$ the derivative of $\langle\langle X,Y\rangle\rangle$ along solutions of $Z$ is denoted by $\mathcal{L}_Z\langle\langle X,Y\rangle\rangle$.
For a Riemannian or singular Riemannian metric one can show that there exists a unique 1-form field, $\mathbb{I}\nabla_XY\in T^*G$ that satisfies the following properties \cite{Stocia}:
\begin{align}
\mathcal{L}_Z\langle\langle X,Y\rangle\rangle&=\langle \mathbb{I}\nabla_ZX,Y\rangle+\langle \mathbb{I}\nabla_ZY,X\rangle,\label{eq:Metric}\\
\mathbb{I}\nabla_XY-\mathbb{I}\nabla_YX&=\mathbb{I}[X,Y].\label{eq:Symmetric}
\end{align}
For a given $X,Y\in TG$ this 1-form field is explicitly given by the Koszul formula:
{\small
\begin{align*}
\langle\mathbb{I}\nabla_XY,Z\rangle&=\frac{1}{2}\left(\mathcal{L}_X\langle\langle Y, Z\rangle\rangle+\mathcal{L}_Y\langle\langle Z, X\rangle\rangle-\mathcal{L}_Z\langle\langle X,Y\rangle\rangle
\right.\\
&\:\:\:\left.-\langle\langle X,[Y, Z]\rangle\rangle+\langle\langle Y, [Z,X]\rangle\rangle+\langle\langle Z, [X,Y]\rangle\rangle\right).
\end{align*}
}
It can be shown that this allows one to define a covariant derivative called the \textit{lower derivative} \cite{Stocia}, $\mathbb{I}\nabla : TG\times TG \mapsto T^*G$, that takes values in $T^*G$.

If the metric is Riemannian then $\mathbb{I}$ is an isomorphism and then $\nabla_X Y\triangleq \mathbb{I}^{-1}(\mathbb{I}\nabla_XY)$ defines a unique \emph{covariant derivative} or \emph{connection},  called the Levi-Civita connection. In this case (\ref{eq:Metric}) states that the connection is metric and (\ref{eq:Symmetric}) states that the connection is symmetric or torsion free.
From a mechanical system point of view what is more crucial is the property of metricity given by (\ref{eq:Metric}).
Using these notations one can write down Newtons equations in the intrinsic fashion 
\begin{align}
\mathbb{I}\nabla_{\dot{g}}\dot{g}&=\gamma_g,\label{eq:MechanicalSystem}
\end{align} 
where $\gamma_g$ is the generalized force acting on the system. This expression is valid for both Riemannian metrics as well as singular Riemannian metrics. Since $\mathbb{I}$ is an isomorphism for  Riemannian (non-singular) metrics the mechanical system can be alternatively written as $\nabla_{\dot{g}}\dot{g}=\mathbb{I}^{-1}\gamma_g=\Gamma_g$.  Where as for singular Riemannian metrics the corresponding mechanical system can not be written in this fashion. In the Riemannian case $\nabla_{\dot{g}}\dot{g}$ has the notion of intrinsic acceleration.

%%%%%%%%%%%%%%%%%%%%%%%%%%%%%%%%%%%%%%%%%
%%%%%%%%%%%%%%%%%%%%%%%%%%%%%%%%%%%%%%%%%
\subsubsection{Mechanical Systems on the Circle}\label{Secn:MechSysOnS}
In this section we consider the special class of mechanical systems that evolve on the Lie-group $\mathbb{S}$ with the desired output evolving on another Lie-group $G$. 

The space of all possible vector fields on $\mathbb{S}$, referred to as the tangent bundle to $\mathbb{S}$ is denoted by $T\,\mathbb{S}\equiv \mathbb{S}\times \mathbb{R}$. The kinetic energy of the system defines a Riemannian metric on $\mathbb{S}$. It  is defined by $\langle\langle \zeta,\eta\rangle\rangle=\mathbb{I}(\theta)\zeta\eta$ where the inertia $\mathbb{I}(\theta)>0$ and $\zeta,\eta\in \mathbb{R}$. The unique Levi-Civita connection that corresponds to the Riemannian metric
$\langle\langle \zeta,\eta\rangle\rangle$ on $\mathbb{S}$ is explicitly given by
\begin{align*}
\nabla_\zeta \eta &=d\eta(\zeta)+\Gamma^1_{11}(\theta)\, \zeta \eta,
\end{align*}
where
\begin{align*}
\Gamma^1_{11}(\theta)&=\frac{1}{2\mathbb{I}}\dfrac{\partial \,\mathbb{I}}{\partial \theta}.
\end{align*}
The significance of the Levi-Civita connection is that it satisfies the metricity condition given by
\begin{align*}
\mathcal{L}_\xi\langle\langle \zeta,\eta\rangle\rangle&=\langle\mathbb{I}\nabla_\xi \zeta,\eta\rangle+\langle\mathbb{I}\nabla_\xi \eta,\zeta\rangle,
\end{align*}
for vector fields $\xi(\theta),\zeta(\theta),\eta(\theta)\in T\mathbb{S}$.
A mechanical system on $\mathbb{S}$ with kinetic energy equal to $\frac{1}{2}\langle\langle \omega,\omega\rangle\rangle=\frac{1}{2}\mathbb{I}(\theta)\omega^2$ is then intrinsically represented by
\begin{align}
\dot{\theta}&=\omega,\label{eq:KinematicsS}\\
\mathbb{I}\,\nabla_\omega\omega &=\tau\label{eq:MechanicalSystemS}
\end{align}
where $\tau\in \mathcal{R}$ is the generalized force. Explicitly we have that (\ref{eq:MechanicalSystemS}) is given by 
$\mathbb{I}\,\dot{\omega}+ \mathbb{I}\,\Gamma^1_{11}(\theta)\, \omega^2=\tau$.

%%%%%%%%%%%%%%%%%%%%%%%%%%%%%%%%%%%%%%%%%

\begin{proofoftheorem}
Consider any small $\epsilon>0$ and compact subsets $\mathcal{X}_s\subset G_y\times \mathcal{G}_s\times \mathcal{G}_s$ and $\mathcal{X}_a\subset  G_a\times \mathcal{G}_a$. 

Let $V_{y} : G_y \mapsto \mathbb{R}$ be a polar Morse function on $G_y$ with the unique minimum at $y_0$. 
Let ${\eta_s}$ be the gradient of $V_s\triangleq V_y\circ y: G_s \mapsto \mathbb{R}$.  
Consider the function $W_s : G_y\times \mathcal{G}_s\times \mathcal{G}_s\times  \mathcal{G}_a$ explicitly given by
\begin{multline}
W_s=k_pV_{s}(g_s)+\frac{1}{2}\langle\langle{v_s}, {v_s}\rangle\rangle_s+\frac{\gamma}{2}\langle\langle {v_I}, {v_I}\rangle\rangle_s\\
%\nonumber\\&\:\:\:\:
+\alpha\langle\langle{\eta_s}, {v_s}\rangle\rangle_s+\beta \langle\langle {v_I}, {v_s}\rangle\rangle_s+\sigma \langle\langle {v_I}, {\eta_s}\rangle\rangle_s.\label{eq:W_s}
\end{multline}

Let $z_s=[||{v_I}|| \:\:\:\: ||{\eta_s}|| \:\:\:\: ||{v_s}||]$. Let $\mathcal{W}_u$ be the set contained by the smallest level set of $W_s$ containing  
$\mathcal{X}_s$ and let $\mathcal{W}_l$ be the set contained by the largest level set of $W_s$ contained in the set where $||z_s||<\epsilon$. Let $k_{s}>0$ be the smallest values such that $||z_s||\leq k_{s}$ for all $(y,{v_s},{v_I})\in \mathcal{W}_u$ and $||v_a||\leq k_{a}$ for all ${v_a}\in \mathcal{X}_a$.

Let $\vartheta$ be such that $\langle\langle\eta_s,\eta_s\rangle\rangle/(2\vartheta) \leq V_y(y)$ on $\mathcal{W}_u$. 
The existence of such a $\vartheta$ is guaranteed by the assumption that $V_{y}$ is a polar Morse function. 
Then it follows that
$W\geq \frac{1}{2} z_s^T P_s z_s\geq \frac{1}{2} \lambda_{\min}(P_s) ||z_s||^2$ where 
\begin{align*}
P_s=\begin{bmatrix}
\gamma &  -\sigma & -\beta\\
-\sigma & \frac{k_p}{\vartheta}  & -\alpha \\
-\beta  &-\alpha & 1
\end{bmatrix}.
\end{align*}
Let $ \mu_{\mathrm{min}}<||\mathbb{I}\nabla\eta|| < \mu_{\mathrm{max}}$ for some $\mu_{\mathrm{min}},\mu_{\mathrm{max}}>0$ and $\mu_a=\mathrm{min}\,||\mathbb{I}_a\nabla\eta_a||$ on $\mathcal{W}_u$. Since $V_y$ is assumed to be a polar Morse function these bounds are guaranteed to exist.
When one chooses $\beta=\frac{k_I}{k_d}$, 
$\sigma=2\kappa {k_I}$, $\gamma=\frac{k_I(\alpha k_d+k_p)}{k_d}$, then it can be shown that if the controller gains are chosen to satisfy (\ref{eq:kICond})--(\ref{eq:kpCond})
then $P_s$ is positive definite.

Differentiating $W_s$ along the dynamics of the closed loop system we have
\begin{align*}
\dot{W}_s
&\leq-z_s^TQ_sz_s\\&\:\:\:\:
- \langle\epsilon_p (k_p\eta_s+k_d\omega_e+k_Iv_I),v_s+\alpha\eta_s+\beta v_I\rangle\\
&\:\:\:\:-\langle B\mathbb{I} _s(k_p\eta_s+k_dv_s+k_Iv_I),v_a\rangle\\
&\:\:\:\:+ \langle \tau_P^s,{v_s}+\alpha{\eta_s}+\beta {v_I}\rangle
+ \langle \tau_s+\Delta_s,{v_s}+\alpha{\eta_s}+\beta {v_I}\rangle,
\end{align*}
where for $\delta \triangleq (1-\mu_{\min}/\mu_{\max})$ 
 \begin{align*}
 Q_s&=\begin{bmatrix}
 \frac{k_I^2}{k_d} & 0 & -\delta k_I\\
 0 &   \left(\alpha k_p-\frac{2 k_d}{\mu_{\max}}\right) & \frac{(k_I-\alpha k_d^2)}{2k_d} \\
-\delta k_I & \frac{(k_I-\alpha k_d^2)}{2k_d} & k_d- {\alpha\mu_{\max}}
 \end{bmatrix}.
 \end{align*}
It can be shown that we can pick gains such that $\lambda_{\min}(Q_s)$ is arbitrary.

The assumption that $\tau_s$ is a quadratic velocity terms implies that the term $\langle \tau_s,{v_s}+\alpha{\eta_s}+\beta {v_I}\rangle$ is cubic in the velocities.
Thus we see that there exists $g_1,g_2,g_3\geq 0$ such that
\[
\langle \tau_s,{v_s}+\alpha{\eta_s}+\beta {v_I}\rangle\leq g_1||z_s||^3+g_2||v_a||||z_s||^2+ g_3||v_a||^2||z_s||.
\]
Also let $g_{\tau}^s>0$ be such that $\tau^s_P(g_s,g_a)<g_{\tau}^s$.
Let $k_a^0>0$ be such that ${(2V_a(g_a)+||v_a||^2)}<k_a^0$ on $\mathcal{X}_a$ and $k_a>k_a^0$. We will show that it is possible to ensure  
$\lim_{t\to \infty} (y(t),{v_s}(t))=(\bar{y}_0,0)$ semi-globally and locally exponentially while $||v_a(t)||<k_a$ for all $t>0$.

For all $(y,v_a,v_I)\in \mathcal{X}_u$ and $v_a$ such that $||v_a||\leq k_a$ we have 
\begin{multline}
\dot{W}_s \leq -(\lambda_{\mathrm{min}}(Q_s)-g_0||\epsilon_{p}||-g_1k_s-g_2k_{a})||z_s||^2\\+(3||\Delta_s||+g_\tau^s+g_3k_{a}^2)\,||z_s||,\nonumber
\end{multline}
where   $g_0\triangleq \mathrm{max}\{k_p,k_d\}$, $g_1\triangleq (3||\Delta_s||+g_\tau^s+g_3k_a^2)$,
and $\epsilon_{p}$ is an operator that depends on the magnitude of the parameter uncertainties of $\mathbb{I}_s$ and $B$.
The right hand side is less than zero when
\[
||z_s||>\epsilon_c=\frac{(3||\Delta_s||+g_\tau^s+g_3k_{a}^2)}{(\lambda_{\mathrm{min}}(Q_s)-g_0||\epsilon_{p}||-g_1k_s-g_2k_{a})}.
\]
Recall that we have shown that the gains $k_p,k_I,k_d$ can be chosen sufficiently large so that $\lambda_{\mathrm{min}}(Q_s)$ can be made arbitrarily large. Thus it can be shown that if the parametric uncertainty is small so that $||\epsilon_p||<1$ then the gains $k_p,k_I,k_d$ can be chosen sufficiently large so that $\epsilon_c<\epsilon$ for any given $k_{s},\epsilon>0$.
Therefore we have shown that that there exists gains $k_p,k_I,k_d$ such that $W_s$ is strictly decreasing in $\mathcal{W}_u/\mathcal{W}_l$ and hence, provided that $||v_a(t)||\leq k_a$ for all time $t>0$, 
the trajectories of the closed loop system can be made to converge to the set $\mathcal{W}_l$. Since $W_s$ is quadratically bounded from below it follows that the convergence is exponential.

If the disturbances and the uncertainties are constant the Lasalle's invariance principle implies that the trajectories of (\ref{eq:System_s}) converge to the largest invariant set contained in the set where $\dot{W}_s\equiv 0$ contained in $\mathcal{W}_l$. For mechanical systems with constant unknown disturbances and constant uncertainty these invariant sets are of the form $(\bar{y},0,\bar{v}_I)$ where $\bar{y}$ is a critical point of $V_y$, and $\bar{v}_I$ is a constant. 
Thus proving that $\lim_{t\to \infty} (y(t),{v_s}(t))=(\bar{y}_0,0)$ semi-globally and locally exponentially in the presence of bounded parametric uncertainty and bounded constant disturbances for all initial conditions in $\mathcal{X}_s$ other than the unstable equilibria and their stable manifolds provided that $||v_a(t)||\leq k_a$ for all $t>0$.

The exponential convergence implies that there exists $\kappa>0$ such that $||v_s(t)||\leq ||v_s(0)||e^{-\kappa t}$ and $||\tau_u-\bar{\tau}_u||\leq \nu\, e^{-\kappa t}$ where $\nu=(||\Delta_s||+g_{\tau}^s+3k_sg_0)$. 
The assumption that $\tau_a(v_s,v_a)$ is quadratic in the velocities implies that there exists $g^a_1\geq 0$ such that $\langle  \tau_{a}, {v_a}\rangle\leq g_3||v_s||||v_a||^2$.
Consider the non-negative function $W_a: G_a\times \mathcal{G}_a\mapsto \mathbb{R}$ defined to be $W_a=V_a+{\langle \langle v_a,v_a\rangle\rangle}/2$.
The derivative of this function along the dynamics of the closed loop system satisfies
$\dot{W}_a\leq 2\left(\nu||B|| +k_sg^a_1\right)e^{-\kappa t}\,W_a$.
This gives that 
\begin{align*}
W_a&\leq W_a(0)\,\exp{\left(\frac{2\left(\nu||B||+k_sg^a_1\right)}{\kappa}\right)}.
\end{align*}
Hence we have that $||v_a(t)||\leq {(2V_a(0)+||v_a(0)||^2)}\,\exp{\left(\frac{\left(\nu||B||+k_sg^a_1\right)}{\kappa}\right)}$. Let $k_a^0>0$ be such that ${(2V_a(g_a)+||v_a||^2)}\leq k_a^0$ on $\mathcal{X}_a$.
We have shown that by picking sufficiently large PID gains $k_p,k_I,k_d$ one can make $\lambda_{\min}(Q_s)$ and hence $\kappa$ sufficiently large. Thus there exists gains such that 
${(2V_a(0)+||v_a(0)||^2)}\,e^{\frac{\left(\nu||B||+k_sg^a_1\right)}{\kappa}}$ can be made less than $k_a>k_a^0$ for all $(g_a,v_a)\in \mathcal{X}_a$ and hence ensure that $||v_a(t)||<k_a$ for all time $t>0$.

\end{proofoftheorem}

%%%%%

%%%%%%%%%%%%%%%%%%%%%%%%%%%%%%%%%%%%%%%%%%%%%%

\end{appendix}

\bibliographystyle{IEEEtran}

\bibliography{RollingHoop}
\newpage
\end{document}